\documentclass[12pt]{article}

\usepackage{amssymb}
\usepackage{latexsym}
\usepackage{exscale}

\title{Vector-valued Hausdorff-Young inequality \\
on compact groups}

\author{{\sc J. Garc\'{\i}a-Cuerva} and {\sc J.Parcet}}

\date{}

\newtheorem{The}{Theorem}[section]
\newtheorem{Lem}[The]{Lemma}
\newtheorem{Pro}[The]{Proposition}
\newtheorem{Cor}[The]{Corollary}
\newtheorem{Def}[The]{Definition}
\newtheorem{Rem}[The]{Remark}

\newcommand{\fin}{\hspace*{\fill} $\Box$ \vskip0.2cm}

\newcommand{\R}{\mathbb{R}}
\newcommand{\C}{\mathbb{C}}
\newcommand{\G}{\widehat{G}}
\newcommand{\B}{\mathcal{B}}
\newcommand{\Hi}{\mathcal{H}}
\newcommand{\Le}{\mathcal{L}}
\newcommand{\Co}{\mathcal{C}}
\newcommand{\bcheck}[1]{{\mbox{\Large$\check{\raisebox{.5em}{}}$}}\!\!\!#1}

\begin{document}

\textwidth=13.5cm

\maketitle

 \footnote{Research supported in part by the TMR Network
`Harmonic Analysis'.} \footnote {Research supported in part by
Project BFM 2001/0189, Spain.}


\section*{Introduction}
\label{section1}

\large\normalsize

If $1 \le p \le 2$ and $p' = p / (p-1)$ denotes the conjugate
exponent of $p$, the classical Hausdorff-Young inequality
establishes the boundedness of the Fourier transform from
$L^p(\R)$ into $L^{p'}(\R)$. Its proof is obtained by complex
interpolation between the obvious case $p = 1$ and the case $p =
2$ given by Plancherel theorem. In the same spirit, Kunze applied
in $1958$ new techniques of non-commutative integration introduced
by Dixmier \cite{D} and Segal \cite{S1,S2} to study this
inequality on locally compact unimodular groups, see \cite{Ku}. In
particular, for a compact non necessarily abelian group $G$, he
proved the boundedness of the Fourier transform from $L^p(G)$ to
$\mathcal{L}^{p'}(\G)$. That is, $$\Big( \sum_{\pi \in \G} d_{\pi}
\|\widehat{f}(\pi)\|_{S_{d_{\pi}}^{p'}}^{p'} \Big)^{1/p'} \le
\Big( \int_G |f(g)|^p d \mu(g) \Big)^{1/p} \ \ \ \mbox{for} \ \ 1
< p \le 2,$$ and with the obvious modifications for $p = 1$. Here
$\pi \in \G$ denotes an irreducible unitary representation of $G$,
$d_{\pi}$ is the degree of $\pi$, $S_n^{p'}$ stands for the $n
\times n$ dimensional Schatten class of exponent $p'$ and $\mu$ is
the Haar measure of $G$ normalized so that $\mu(G) = 1$. On the
other hand, Peetre presented in 1969 the first work \cite{Pe}
analyzing the Hausdorff-Young inequality for Banach valued
functions $f: \R \rightarrow B$. In this case, the validity of the
inequality for some fixed $p$ depends on the Banach space $B$.
This leads to the notion of Fourier type of a Banach space with
respect to a locally compact abelian group, introduced by Milman
in \cite{M}. The theory of Fourier type with respect to locally
compact abelian groups was further developed by several authors,
see \cite{A}, \cite{B1}, \cite{B2}, \cite{GKKT}, \cite{GKK} and
\cite{Ko}.

However, as far as we know, in the non-commutative setting there
is no analogous theory of Fourier type described in the
literature. Our aim is to fill this gap. Namely, to analyze the
validity of Kunze's results for vector-valued functions. In this
work we investigate the validity of the Hausdorff-Young inequality
for vector-valued functions defined on a compact group. As can be
seen throughout the paper, compactness is an essential assumption
in many of the results we present here. For a non-commutative
compact group $G$, the vector-valued Fourier transform must be
defined for irreducible unitary representations $\pi \in \G$ and
its values are vector-valued matrices. Therefore, just to start
talking about the Hausdorff-Young inequality, one has to be able
to define norms for vector-valued matrices. By Ruan's theorem
\cite{R}, this matricial structure leads us to consider an
operator space structure on the vector space where we are taking
values. It appears clear that, in order to develop a theory of
Fourier type in this context, we shall need to take values in
operator spaces rather than Banach spaces. This crucial point is
obviously at the root of the notion of Fourier type.

To conclude, we would like to point out that the theory initiated
in this paper has been further developed in \cite{GMP} and
\cite{GP}. Roughly speaking, the paper \cite{GMP} deals with the
sharpness of Theorems \ref{BochnerLebesgue} and
\ref{PisierSchatten} for compact semisimple Lie groups, see
section \ref{section7} for more on this topic. On the other hand,
the notions of Fourier type and cotype of an operator space with
respect to a compact group are extended in \cite{GP} to the more
general setting of type and cotype with respect to a \lq quantized
orthonormal system\rq. This contains, for instance, the
non-commutative versions of Rademacher or Gauss type and cotype.
All this is used in \cite{GP} to obtain an operator space version
of Kwapie\'n theorem \cite{Kw} characterizing Hilbert spaces by
means of vector-valued orthogonal series.

The organization of the paper is as follows. In section
$\ref{section2}$ we recall the notions of operator space theory
and vector-valued Schatten classes that will be used in the
sequel. In section $\ref{section3}$ we define the Fourier
transform on compact groups for vector-valued functions. We also
study the spaces $\Le_E^p(\G)$, where the Fourier transform takes
values. Some results for which we have not found any reference
have been proved for completeness. In particular, Proposition
\ref{L-inf} and Corollary \ref{L-one} are specially relevant since
they show that the spaces $\Le_E^{\infty}(\G)$ and $\Le_E^1(\G)$
behave with respect to the minimal and projective tensor products
as the classical $L_E^{\infty}$ and $L_E^1$ do with respect to
Grothendieck's tensor norms. Sections $\ref{section4}$ and
$\ref{section5}$ are mainly devoted to showing that the notions of
Fourier type and cotype are well defined, and also to prove some
basic properties. Section $\ref{section5}$ is specially far from
the commutative theory since, as we show there, when dealing with
abelian groups the notion of Fourier cotype reduces to the notion
of Fourier type with respect to the dual group. In section
$\ref{section6}$, given an operator space $E$, we investigate the
Fourier type and cotype of some general operator spaces related to
$E$ such as subspaces, duals, interpolated spaces, etc... Finally,
in section $\ref{section7}$ we investigate the main examples, that
is, Lebesgue spaces and Schatten classes. In particular, for the
vector-valued ones we prove some quantized Minkowski inequalities
that we shall need.

\

\noindent \textbf{Acknowledgment.} We thank Gilles Pisier for some
useful comments.

\section{Operator spaces and Schatten classes}
\label{section2}

The basic theories behind this paper are the theory of operator
spaces and the subsequent theory of vector-valued Schatten
classes. The reader is referred to \cite{ER2} and \cite{P2} for a
basic background on these topics and their connection with the
present work. We begin with a brief summary of the results of
operator space theory that will be used in the sequel.
\begin{itemize}
\item[(a)] \textbf{Definition of operator space}. We will denote
by $\B(\Hi)$ the space of bounded linear operators on some Hilbert
space $\Hi$. For our purposes an \emph{operator space} $E$ can be
defined as a closed subspace of $\B(\Hi)$. Given an operator space
$E \subset \B(\Hi)$ we write $M_n(E)$ for the space $M_n \otimes
E$ of $n \times n$ matrices with entries in $E$ and with the norm
imposed by the natural embedding of $M_n \otimes E$ into
$\B(l_{\Hi}^2(n))$. Here $l_{\Hi}^2(n)$ denotes the Hilbert space
of all $\Hi$-valued $n$-tuples with its natural inner product. On
the other hand, given a vector space $E$ and a collection of norms
$\|-\|_n$ on the spaces $M_n \otimes E$, one can impose some extra
conditions to obtain what is called an \emph{operator space matrix
norm} or \emph{operator space structure} on $E$, see Chapter $2$
of \cite{ER2}. One of the main results of the theory is the
\emph{abstract characterization of operator spaces} given by Ruan
in \cite{R}. Ruan's theorem can be rephrased by saying that for
any operator space structure on a vector space $E$, there exists a
unique Hilbert space $\Hi$ such that the norm from the operator
space structure on $M_n \otimes E$ coincides with the norm induced
by the space $\B(l_{\Hi}^2(n))$.

\item[(b)] \textbf{Complete boundedness}. A linear mapping between
operator spaces $u: E_1 \rightarrow E_2$ is said to be
\emph{completely bounded} if the family of maps $I_{M_n} \otimes
u: M_n(E_1) \rightarrow M_n(E_2)$ satisfy $$\|u\|_{cb} = \sup_{n
\ge 1} \|I_{M_n} \otimes u\|_{\B(M_n(E_1),M_n(E_2))} < \infty.$$
We write $cb(E_1,E_2)$ for the Banach space of completely bounded
maps from $E_1$ to $E_2$ with the $cb$ norm. Let $u \in
cb(E_1,E_2)$, we say that $u$ is a \emph{complete isometry} if the
mappings $I_{M_n} \otimes u$ are isometries for all $n \ge 1$.
Similarly $u$ is called \emph{completely contractive} if
$\|u\|_{cb} \le 1$. We also say that $u$ is a \emph{complete
isomorphism} if it is a completely bounded linear isomorphism
whose inverse is also completely bounded. Finally $u$ is a
\emph{completely isometric isomorphism} if it is also a complete
isometry.

\item[(c)] \textbf{Duality}. Ruan's theorem was used by Blecher
and Paulsen in \cite{BP} and by Effros and Ruan in \cite{ER1} to
get a \emph{duality theory} in the category of operator spaces. It
was shown that, by imposing on $M_n \otimes cb(E_1,E_2)$ the norm
induced by $cb(E_1,M_n(E_2))$, we obtain an operator space
structure on $cb(E_1,E_2)$. In particular we have an operator
space structure on the dual space $E^{\star} = cb(E,\C)$. This
notion of duality behaves as Banach space duality in many senses.
For instance, it can be proved that the natural isometric
inclusion $E \subset E^{\star \star}$ is a complete isometry.

\item[(d)] \textbf{Tensor products}. We are interested in two
tensor norms that will be used repeatedly in this paper. Given two
operator spaces $E_1 \subset \B(\Hi_1)$ and $E_2 \subset
\B(\Hi_2)$ we define their \emph{minimal tensor product} $E_1
\otimes_{\mbox{{\tiny min}}} E_2$ by the natural embedding of $E_1
\otimes E_2$ into $\B(\Hi_1 \otimes_2 \Hi_2)$, where $\otimes_2$
stands for the Hilbertian tensor product. The minimal tensor
product plays the role of the \emph{injective tensor product} of
Banach spaces in the category of operator spaces. Similarly, there
exists an analog for operator spaces of the \emph{projective
tensor product}. It is denoted by $E_1 \otimes^{\wedge} E_2$ and
it was introduced in \cite{BP} and \cite{ER1} independently. The
tensor products $\otimes_{\mbox{{\tiny min}}}$ and
$\otimes^{\wedge}$ are associative and commutative. Here are some
other properties which we will use in the sequel with no further
reference $$\begin{array}{lclcl} E_1 \otimes^{\wedge} E_2 &
\rightarrow & E_1 \otimes_{\mbox{{\tiny min}}} E_2 & \quad &
\mbox{is a comp. contraction}
\\ E_1^{\star} \otimes_{\mbox{{\tiny min}}} E_2 &
\rightarrow & cb(E_1,E_2) & \quad & \mbox{is a comp. isometry}
\\ E_1 \otimes_{\mbox{{\tiny min}}} E_2 & \rightarrow &
cb(E_1^{\star},E_2) & \quad & \mbox{is a comp. isometry}
\\ (E_1 \otimes^{\wedge} E_2)^{\star} & \rightarrow &
cb(E_1,E_2^{\star}) & \quad & \mbox{is a comp. isometric
isomorph}.
\end{array}$$

\item[(e)] \textbf{Complex interpolation}. Let $\{E_0,E_1\}$ be a
compatible couple of Banach spaces in the sense of complex
interpolation. Let us suppose that $E_0$ and $E_1$ have an
operator space structure. In \cite{P1} Pisier showed that, if
$E_{\theta}$ denotes the interpolation space $[E_0,E_1]_{\theta}$,
one can define an operator space structure on $E_{\theta}$ by
imposing on $M_n \otimes E_{\theta}$ the norm of the Banach space
$[M_n(E_0),M_n(E_1)]_{\theta}$. He also proved the analog for
operator spaces of the classical interpolation result for Banach
spaces. Namely, if we assume that $u: E_0 + E_1 \rightarrow F_0 +
F_1$ satisfies the inequalities $\|u\|_{cb(E_0,F_0)} \le C_0$ and
$\|u\|_{cb(E_1,F_1)} \le C_1$, then for $0 < \theta < 1$ we have
the estimate $$\|u\|_{cb(E_{\theta},F_{\theta})} \le
C_0^{1-\theta} C_1^{\theta}.$$
\end{itemize}

We now recall the definition and the main properties of the
Schatten classes. The non-commutative analog of the
$n$-dimensional Lebesgue space $l^p(n)$ is the Schatten class
$S_n^p$ which is defined as the space $M_n$ of $n \times n$
complex matrices with the norm given by
\begin{itemize}
\item[(a)] $\|A\|_{S_n^p} = (\mbox{tr}
|A|^p)^{1/p}$, \ \ if $1 \le p < \infty$.
\item[(b)] $\|A\|_{S_n^{\infty}} = \displaystyle
\sup \Big\{ \|Ax\|_{l^2(n)}: \ \|x\|_{\ell^2(n)} \le 1 \Big\}$, \
\ if $p = \infty$.
\end{itemize}

Now we present the \emph{vector valued Schatten classes}
$S_n^p(E)$, introduced by Pisier in \cite{P2}. The point here is
that the space $E$ where we take values has to be an operator
space. If $p = \infty$, we have by definition $S_n^{\infty} =
\B(l^2(n))$ and so we obtain a natural operator space structure
for $S_n^{\infty}$. We define $S_n^{\infty}(E)$ as the operator
space $S_n^{\infty} \otimes_{\mbox{{\tiny min}}} E$. It is obvious
that $M_n(E)$ and $S_n^{\infty}(E)$ coincide, in what follows we
shall write $S_n^{\infty}(E)$ for $M_n(E)$. If $p = 1$, the
duality $S_n^1 = (S_n^{\infty})^{\star}$ gives a natural operator
space structure on $S_n^1$. We set $S_n^1(E) = S_n^1
\otimes^{\wedge} E$. Finally, since the identity mapping $S_n^1(E)
\rightarrow S_n^{\infty}(E)$ is contractive, we define the classes
$S_n^p(E)$ by means of complex interpolation. Namely, $S_n^p(E) =
[S_n^{\infty}(E),S_n^1(E)]_{1/p}$. The next theorem summarizes
some properties of the vector-valued Schatten classes that will be
used repeatedly throughout the paper, see Chapter $1$ of
\cite{P2}.

\begin{Rem}
\emph{In the same fashion, it is possible to define a natural
operator space structure on the Bochner-Lebesgue spaces, see
Chapter $2$ of \cite{P2}.}
\end{Rem}

\begin{The} [Pisier] \label{Schatten} The vector-valued Schatten
classes satisfy the fo\-llowing properties.
\begin{itemize}
\item[$1.$ ] \textbf{The $cb$ norm}. Let $E_1$ and $E_2$ be operator
spaces and let $1 \le p \le \infty$. Then the $cb$ norm of any
linear mapping $u: E_1 \rightarrow E_2$ is given by $$\|u\|_{cb} =
\sup_{n \ge 1} \|I_{M_n} \otimes
u\|_{\B(S_n^p(E_1),S_n^p(E_2))}.$$
\item[$2.$ ] \textbf{Duality}. Let $1 \le p \le \infty$ and let $p'$
denote the conjugate exponent of $p$. The map $A \in
S_n^{p'}(E^{\star}) \mapsto \textnormal{tr} (A \hspace{2pt} \cdot)
\in S_n^p(E)^{\star}$ is completely isometric.
\item[$3.$ ] \textbf{Complex interpolation}. Let $1 \le p_0,p_1
\le \infty$, $0 < \theta < 1$ and assume that $\{E_0,E_1\}$ is a
compatible couple of operator spaces. Then, letting
$p_{\theta}^{-1} = (1 - \theta) p_0^{-1} + \theta p_1^{-1}$, we
have $$[S_n^{p_0}(E_0),S_n^{p_1}(E_1)]_{\theta} =
S_n^{p_{\theta}}(E_{\theta}).$$
\item[$4.$ ] \textbf{Ordered norms}. Let $1 \le p_1
\le p_2 \le \infty$. Then the identity map $S_n^{p_1}(E)
\rightarrow S_n^{p_2}(E)$ is a contraction.
\item[$5.$ ] \textbf{Fubini type theorems}. Let $1 \le p \le \infty$
and let $n_1, n_2, n \ge 1$. Then, we have completely
isometrically $$S_{n_1}^p(S_{n_2}^p(E)) \simeq
S_{n_2}^p(S_{n_1}^p(E)) \quad \mbox{and} \quad
S_n^p(L_E^p(\Omega)) \simeq L_{S_n^p(E)}^p(\Omega).$$
\end{itemize}
\end{The}

\section{Vector-valued Fourier transform}
\label{section3}

We assume the reader is familiar with the language of
non-commutative abstract harmonic analysis on compact groups. In
any case all the results we use here can be found in \cite{F}. In
what follows we shall assume that $G$ is a compact Hausdorff
topological group endowed with its Haar measure $\mu$ normalized
so that $\mu(G) = 1$. The mapping $\pi: G \rightarrow
U(\C^{d_{\pi}})$ will denote an irreducible unitary representation
of $G$ of degree $d_{\pi}$. That is, $\pi \in \G$ where the symbol
$\G$ stands for the dual object of $G$.

\begin{Def} \emph{Given an operator space $E$, $f \in L_E^1(G)$
and $\pi \in \G$, the vector-valued \textbf{Fourier coefficient}
of $f$ at $\pi$ is defined as the operator}
$$\widehat{f}(\pi) = \int_G f(g) \pi(g)^{\star} d\mu(g) \in
\B(\C^{d_{\pi}},E^{d_{\pi}}).$$
\end{Def}

We interpret this operator-valued integral in the weak sense. That
is, given an orthonormal basis $\{v_1,v_2 , \ldots, v_{d_{\pi}}\}$
of $\C^{d_{\pi}}$ and $u \in \C^{d_{\pi}}$, we define the $j$-th
component of $\widehat{f}(\pi) u$, with respect to that basis, by
the following element of $E$ $$\int_G f(g) \langle\pi(g)^{\star}
u, v_j\rangle d \mu(g).$$ Since $\pi(g)$ is unitary, it follows
that the vector-valued Fourier coefficients are well-defined for
all $f$ in $L_E^1(G)$. Once we have fixed the basis of
$\C^{d_{\pi}}$, we can identify $\B(\C^{d_{\pi}},E^{d_{\pi}})$
with the space $M_{d_{\pi}} \otimes E$. This leads us to write the
\textbf{Fourier transform} operator $\mathcal{F}_{G,E}$, in the
form $$\mathcal{F}_{G,E}: L_E^1(G) \longrightarrow \prod_{\pi \in
\G} M_{d_{\pi}} \otimes E.$$

The first step to study the Hausdorff-Young inequality is to find
a \emph{natural} $L^p$ norm for this Cartesian product, which we
denote by $\mathcal{M}_E(\G)$.

\begin{Def} \emph{Let $E$ be an operator space and $1 \le p <
\infty$, the \textbf{spaces} $\Le_E^p(\G)$ are defined as follows}
\begin{eqnarray*}
\Le_E^p(\G) & = & \Big\{ A \in \mathcal{M}_E(\G): \ \
\|A\|_{\Le_E^p(\G)} = \Big( \sum_{\pi \in \G} d_{\pi}
\|A^{\pi}\|_{S_{d_{\pi}}^p(E)}^p \Big)^{1/p} < \infty \Big\}
\\ \Le_E^{\infty}(\G) & = & \Big\{ A \in \mathcal{M}_E(\G): \ \
\|A\|_{\Le_E^{\infty}(\G)} = \sup_{\pi \in \G}
\|A^{\pi}\|_{S_{d_{\pi}}^{\infty}(E)} < \infty \Big\}.
\end{eqnarray*}
\end{Def}

We write $\Le^p(\G)$ for the case $E = \C$. Note that we require
the vector space $E$ to be an operator space. This condition is
necessary since we are making use of the spaces $S_n^p(E)$, see
Pisier's monograph \cite{P2} for more on this topic. The family of
spaces $\Le_E^p(\G)$ is a particular case of a bigger family of
spaces studied in Chapter $2$ of \cite{P2}, the spaces
$\ell_p(\mu,\{E_i\})$. This remark allows us to provide the spaces
$\Le_E^p(\G)$ with the natural operator space structure induced by
$\ell_p(\mu,\{E_i\})$. We now summarize the main properties of
these spaces.

\begin{itemize}
\item[(a)] \textbf{Duality}. Let $1 \le p < \infty$ and let $p'$
be the conjugate exponent of $p$. Then the following map is a
completely isometric isomorphism $$A \in \Le_{E^{\star}}^{p'}(\G)
\longmapsto \sum_{\pi \in \G} d_{\pi} \mbox{tr} (A^{\pi} \cdot)
\in \Le_E^p(\G)^{\star}.$$
\item[(b)] \textbf{Complex interpolation}. Let $1 \le p_0,p_1
\le \infty$. Assume that $\{E_0,E_1\}$ is a compatible couple of
operator spaces. Then
$\{\Le_{E_0}^{p_0}(\G),\Le_{E_1}^{p_1}(\G)\}$ is also a compatible
couple and, for $0 < \theta < 1$ and $p_{\theta}^{-1} = (1 -
\theta) p_0^{-1} + \theta p_1^{-1}$, we have that
$$[\Le_{E_0}^{p_0}(\G),\Le_{E_1}^{p_1}(\G)]_{\theta} =
\Le_{E_{\theta}}^{p_{\theta}}(\G).$$
\item[(c)] \textbf{Ordered norms}. The embedding
$\Le_E^{p_1}(\G) \rightarrow \Le_E^{p_2}(\G)$ is contractive
whenever $1 \le p_1 \le p_2 \le \infty$.
\item[(d)] \textbf{Fubini type theorems}. Let $1 \le p \le \infty$
and let $n \ge 1$. Then the following are completely isometric
isomorphisms $$S_n^p(\Le_E^p(\G)) \simeq \Le_{S_n^p(E)}^p(\G)
\quad \mbox{and} \quad \Le_{L_E^p(\Omega)}^p(\G) \simeq
L_{\Le_E^p(\G)}^p(\Omega).$$
\end{itemize}

We now present a couple of results concerning these spaces for
which we have not found any reference. These will be especially
useful in the study of the Fourier cotype, see for instance the
proof of Proposition \ref{CotypeInfinite}. We first need a
technical result, which is an inequality of H\"{o}lder type.

\begin{Lem} \label{Holder} Let $E$ be an operator space, $n_1, n_2
\ge 1$ and $1 \le p \le \infty$. Let us consider $A \in M_{n_1}
\otimes E$ and $B_{ij} \in M_{n_1}$ for $1 \le i,j \le n_2$. Then
$$\Big\| \Big( \,\ \textnormal{tr} (A B_{ij}) \,\ \Big)
\Big\|_{S_{n_2}^1(E)} \le \|A\|_{S_{n_1}^{p'}(E)} \Big\| \Big( \,\
B_{ij} \,\ \Big) \Big\|_{S_{n_1}^p(S_{n_2}^1)}.$$
\end{Lem}

\emph{Proof}. If $a_{kl}$ and $b_{kl}^{ij}$ denote the entries of
$A$ and $B_{ij}$ respectively, then we can write $$\Big( \,\
\mbox{tr} (A B_{ij}) \,\ \Big) = \sum_{k,l = 1}^{n_1} \Big( \,\
b_{lk}^{ij} \,\ \Big) \otimes a_{kl} \,\ \in \,\ M_{n_2} \otimes
E.$$ Hence, recalling the completely isometric isomorphism from
$cb(S_n^1,E^{\star})$ onto $(S_n^1 \otimes^{\wedge} E)^{\star}$
given by $\Psi(A \otimes e) = \Phi(A)(e)$, we obtain
\begin{eqnarray*}
\Big\| \Big( \,\ \mbox{tr} (A B_{ij}) \,\ \Big)
\Big\|_{S_{n_2}^1(E)} & = &
\sup_{\|\Phi\|_{cb(S_{n_2}^1,E^{\star})} \le 1} \, \Big| \sum_{k,l
= 1}^{n_1} \big[ \Phi \Big( \,\ b_{lk}^{ij} \,\ \Big) \big]
(a_{kl}) \Big| \\ & = & \sup_{\|\Phi\|_{cb(S_{n_2}^1,E^{\star})}
\le 1} \, \big| \mbox{tr} (C A) \big| \\ & \le &
\|A\|_{S_{n_1}^{p'}(E)} \,\
\sup_{\|\Phi\|_{cb(S_{n_2}^1,E^{\star})} \le 1} \ \
\|C\|_{S_{n_1}^p(E^{\star})} \\ & \le & \|A\|_{S_{n_1}^{p'}(E)} \
\ \Big\| \Big( \,\ B_{ij} \,\ \Big) \Big\|_{S_{n_1}^p(S_{n_2}^1)}
\end{eqnarray*}
where $C = \big[ I_{M_{n_1}} \otimes \Phi \big] \Big( \,\ B_{ij}
\,\ \Big) \in M_{n_1} \otimes E^{\star}$. This completes the
proof. \fin

For fixed $\pi_0 \in \G$ and $1 \le i_0,j_0 \le d_{\pi_0}$, we
define $M(\pi_0,i_0,j_0) \in \mathcal{M}_{\C}(\G)$ by the
relations $M(\pi_0,i_0,j_0)_{ij}^{\pi} = \delta_{\pi \pi_0}
\delta_{i i_0} \delta_{j j_0}$.

\begin{Pro} \label{L-inf} The following is a completely isometric
isomorphism $$\begin{array}{crcl} \Lambda: & \Le_E^{\infty}(\G) &
\longrightarrow & cb(\Le^1(\G),E) \\ & A & \longmapsto & \sum_{\pi
\in \G} d_{\pi} \textnormal{tr}(A^{\pi} \cdot) \end{array}$$
\end{Pro}

\emph{Proof}. We just need to show that $\Lambda$ is an isometric
isomorphism, since we have the natural isometric isomorphisms
$$S_n^{\infty}(\Le_E^{\infty}(\G)) \sim
\Le_{S_n^{\infty}(E)}^{\infty}(\G) \quad \mbox{and} \quad
cb(\Le^1(\G),S_n^{\infty}(E)) \sim
S_n^{\infty}(cb(\Le^1(\G),E)).$$

\noindent \textbf{1. $\Lambda$ is a contraction}. By expressing
the $cb$ norm in terms of the Schatten class $S^1$, we have
\begin{eqnarray*}
\lefteqn{\|\Lambda(A)\|_{cb(\Le^1(\G),E)}}
\\ & \le & \sup_{n \ge 1} \Big\{ \sum_{\pi \in \G} d_{\pi} \Big\|
\Big( \,\ \mbox{tr} (A^{\pi} B_{ij}^{\pi}) \,\ \Big)
\Big\|_{S_n^1(E)}: \ \ \Big\| \Big( \,\ B_{ij} \,\ \Big)
\Big\|_{S_n^1(\Le^1(\G))} \le 1 \Big\}
\end{eqnarray*}
But Lemma \ref{Holder} with $p = 1$ gives
\begin{eqnarray*}
\sum_{\pi \in \G} d_{\pi} \Big\| \Big( \,\ \mbox{tr} (A^{\pi}
B_{ij}^{\pi}) \,\ \Big) \Big\|_{S_n^1(E)} & \le & \sum_{\pi \in
\G} d_{\pi} \|A^{\pi}\|_{S_{d_{\pi}}^{\infty}(E)} \Big\| \Big( \,\
B_{ij}^{\pi} \,\ \Big) \Big\|_{S_{d_{\pi}n}^1} \\ & \le &
\|A\|_{\Le_E^{\infty}(\G)} \Big\| \Big( \,\ B_{ij} \,\ \Big)
\Big\|_{S_n^1(\Le^1(\G))}.
\end{eqnarray*}

\noindent \textbf{2. $\Lambda$ is an isometry}. For fixed $\pi \in
\G$, we define $B(\pi,i,j)$ to be the element of $\Le^1(\G)$ given
by $d_{\pi}^{-1} M(\pi,j,i)$ and we denote by $\mathbf{B}(\pi)$
the matrix with entries $B(\pi,i,j)$ where $1\le i,j \le d_{\pi}$.
Due to the natural complete isometry $\Le^1(\G) \rightarrow
cb(\Le^{\infty}(\G),\C)$, it is not difficult to check that
$\|\mathbf{B}(\pi)\|_{S_{d_{\pi}}^{\infty}(\Le^1(\G))} = 1$. Since
this works for any $\pi \in \G$, we get
$$\|\Lambda(A)\|_{cb(\Le^1(\G),E)} \ge \sup_{\pi \in \G}
\|[I_{M_{d_{\pi}}} \otimes \Lambda(A)]
(\mathbf{B}(\pi))\|_{S_{d_{\pi}}^{\infty}(E)} =
\|A\|_{\Le_E^{\infty} (\G)}.$$

\noindent \textbf{3. $\Lambda$ is surjective}. Let $\Phi \in
cb(\Le^1(\G),E)$, then we define $A \in \mathcal{M}_E(\G)$ by the
relation $$A^{\pi} = \frac{1}{d_{\pi}} \Big( \,\ \Phi(M(\pi,j,i))
\,\ \Big), \quad \pi \in \G.$$ The definition of $A$ gives rise to
the following expression $$\Phi(B) = \sum_{\pi \in \G}
\sum_{i,j=1}^{d_{\pi}} b_{ij}^{\pi} \Phi(M(\pi,i,j)) = \sum_{\pi
\in \G} d_{\pi} \mbox{tr} (A^{\pi} B^{\pi})$$ where $b_{ij}^{\pi}$
are the entries of $B^{\pi}$. Therefore it suffices to check that
$A \in \Le_E^{\infty}(\G)$. But following the notation of Step
$2$, we obtain
\begin{eqnarray*}
\|A^{\pi}\|_{S_{d_{\pi}}^{\infty}(E)} & = & \frac{1}{d_{\pi}}
\Big\| \Big( \,\ \Phi(M(\pi,j,i)) \,\ \Big)
\Big\|_{S_{d_{\pi}}^{\infty}(E)} \\ & \le &
\|\Phi\|_{cb(\Le^1(\G),E)}
\|\mathbf{B}(\pi)\|_{S_{d_{\pi}}^{\infty}(\Le^1(\G))} \le
\|\Phi\|_{cb(\Le^1(\G),E)}.
\end{eqnarray*}
Since $\Phi \in cb(\Le^1(\G),E)$, we have a uniform upper bound.
\fin

The space $\Le_E^{\infty}(\G)$ behaves with respect to the minimal
tensor product as $L_E^{\infty}(\Omega)$ does with respect to the
injective tensor product in the category of Banach spaces. Namely,
as a consequence of Proposition \ref{L-inf}, we have that
$\Le^{\infty}(\G) \otimes_{\mbox{{\tiny min}}} E \hookrightarrow
\Le_E^{\infty}(\G)$ is a complete isometry . The space
$\Le_E^1(\G)$ behaves in the same fashion with respect to the
projective tensor product.

\begin{Cor} \label{L-one} The identity $\Le^1(\G)
\otimes^{\wedge} E \rightarrow \Le_E^1(\G)$ is completely
isometric.
\end{Cor}

\emph{Proof}. By Proposition $3.2.2$ of \cite{ER2} it suffices to
check that the adjoint mapping is a complete isometric
isomorphism. But Proposition \ref{L-inf} gives the following chain
$\Le_E^1(\G)^{\star} \simeq \Le_{E^{\star}}^{\infty}(\G) \simeq
cb(\Le^1(\G),E^{\star}) \simeq (\Le^1(\G) \otimes^{\wedge}
E)^{\star}$ of completely isometric isomorphisms. This completes
the proof. \fin

For the sake of completeness we introduce the \textbf{space
$C_0(\G,E)$}. It is defined as the collection of those $A \in
\mathcal{L}_E^{\infty}(\G)$ satisfying $$ \forall \ \ \varepsilon
> 0 \quad \mbox{we have} \quad
\|A^{\pi}\|_{S_{d_{\pi}}^{\infty}(E)} < \varepsilon \quad
\mbox{except for finitely many} \ \ \pi \in \G.$$ As a subspace of
$\Le_E^{\infty}(\G)$ this space inherits a natural operator space
structure. The only two results about the spaces $\Le_E^p(\G)$
that fail at $p = \infty$ are the density of $\Le^p(\G) \otimes E$
in $\Le_E^p(\G)$ and duality, the predual of
$\Le_{E^{\star}}^1(\G)$ is not $\Le_E^{\infty}(\G)$. However, it
is easy to see the density of $C_0(\G) \otimes E$ in $C_0(\G,E)$.
On the other hand, the dual of $C_0(\G,E)$ is completely
isomorphic to $\Le_{E^{\star}}^1(\G)$.

\begin{Pro} The following is a completely isometric isomorphism
$$\begin{array}{crcl} \Lambda: & \Le_{E^{\star}}^1(\G) &
\longrightarrow & C_0(\G,E)^{\star} \\ & A & \longmapsto &
\sum_{\pi \in \G} d_{\pi} \textnormal{tr} (A^{\pi} \cdot)
\end{array}$$
\end{Pro}

\emph{Proof}. Taking into account the natural isometric
isomorphisms given by $S_n^1(\Le_{E^{\star}}^1(\G)) \sim
\Le_{S_n^{\infty}(E)^{\star}}^1(\G)$ and
$C_0(\G,S_n^{\infty}(E))^{\star} \sim S_n^1(C_0(\G,E)^{\star})$ it
is enough to see that $\Lambda$ is an isometric isomorphism. We
prove this fact in several steps.

\vskip5pt

\noindent \textbf{1. $\Lambda$ is a contraction}. This is an
obvious consequence of the duality action on the Schatten classes
$S_n^p(E)$
\begin{eqnarray*}
\|\Lambda(A)\|_{C_0(\G,E)^{\star}} & \le & \sup_{\|B\|_{C_0(\G,E)}
\le 1} \ \ \sum_{\pi \in \G} d_{\pi} |\mbox{tr} (A^{\pi} B^{\pi})|
\\ & \le & \sup_{\|B\|_{C_0(\G,E)} \le 1} \ \ \sum_{\pi \in \G}
d_{\pi} \|A^{\pi}\|_{S_{d_{\pi}}^1(E^{\star})}
\|B^{\pi}\|_{S_{d_{\pi}}^{\infty}(E)} \\ & \le &
\|A\|_{\Le_{E^{\star}}^1(\G)}.
\end{eqnarray*}

\noindent \textbf{2. $\Lambda$ is an isometry}. Let $A \in
\Le_{E^{\star}}^1(\G)$. For all $\varepsilon > 0$ there exists a
finite set $I_{A,\varepsilon} \subset \G$ such that
$$\sum_{\pi \notin I_{A,\varepsilon}} d_{\pi}
\|A^{\pi}\|_{S_{d_{\pi}}^1(E^{\star})} < \varepsilon / 2.$$
Furthermore, for all $\pi \in \G$ there exists
$B_{\varepsilon}^{\pi} \in S_{d_{\pi}}^{\infty}(E)$ of norm $1$
such that $$\mbox{tr} (A^{\pi} B_{\varepsilon}^{\pi}) >
\|A^{\pi}\|_{S_{d_{\pi}}^1(E^{\star})} -
\frac{\varepsilon/2}{|I_{A,\varepsilon}| \displaystyle \max_{\pi
\in I_{A,\varepsilon}} d_{\pi}}$$ where $|I_{A,\varepsilon}|$
denotes the number of elements of $I_{A,\varepsilon}$. Let
$C_{\varepsilon}$ be the element of $C_0(\G,E)$ of norm $1$
defined by $C_{\varepsilon}^{\pi} = B_{\varepsilon}^{\pi}$ if $\pi
\in I_{A, \varepsilon}$ and $0$ otherwise. This Step is completed
by taking $\varepsilon$ arbitrarily small in the expression
\begin{eqnarray*}
\|\Lambda(A)\|_{C_0(\G,E)^{\star}} & \ge & \Big| \sum_{\pi \in \G}
d_{\pi} \mbox{tr} (A^{\pi} C_{\varepsilon}^{\pi}) \Big| >
\|A\|_{\Le_{E^{\star}}^1(\G)} - \varepsilon.
\end{eqnarray*}

\noindent \textbf{3. $\Lambda$ is surjective}. Let $\Phi \in
C_0(\G,E)^{\star}$. Then we define $A \in
\mathcal{M}_{E^{\star}}(\G)$ by the relation $$A^{\pi} =
\frac{1}{d_{\pi}} \Big( \,\ \Phi(M(\pi,j,i) \otimes \,\ \cdot \,\
) \,\ \Big), \quad \pi \in \G.$$ As in Proposition \ref{L-inf}, it
can be shown that $\Phi = \Lambda (A)$ with $A \in
\Le_{E^{\star}}^1(\G)$. \fin

\section{Fourier type}
\label{section4}

Let $E$ be an operator space and let $1 \le p \le 2$. Given $f \in
L^p(G)$ and $e \in E$ it is obvious that the Fourier transform of
$f \otimes e$ coincides with $\widehat{f} \otimes e$. Thus, the
Hausdorff-Young inequality for compact groups (see \cite{Ku} or
Lemma \ref{H-Y} below) provides the relation $\mathcal{F}_{G,E}
(L^p(G) \otimes E) \subset \Le^{p'}(\G) \otimes E$. This motivates
the following definition.

\begin{Def} \emph{Let $1 \le p \le 2$ and let $p'$ be the conjugate
exponent of $p$. We say that the operator space $E$ has
\textbf{Fourier type} $p$ with respect to the compact group $G$ if
the Fourier transform $\mathcal{F}_{G,E}: L^p(G) \otimes E
\rightarrow \Le^{p'}(\G) \otimes E$ can be extended to a
completely bounded operator $$\Lambda_{G,E,p}^1: L_E^p(G)
\longrightarrow \Le_E^{p'}(\G).$$ In that case, we shall denote by
$\Co_p^1(E,G)$ the $cb$ norm of $\Lambda_{G,E,p}^1$.}
\end{Def}

\begin{Rem}
\emph{If the compact group $G$ is also abelian, there exists
already a notion (introduced by Milman in \cite{M}) of Fourier
type of a Banach space with respect to $G$. The only difference
with Milman's notion is that here we require the extended operator
to be completely bounded while in the commutative setting, only
the boundedness of this operator is required.}
\end{Rem}

The first natural question that arises after the definition of
Fourier type is if the extension of $\mathcal{F}_{G,E}$ is always
the natural one. That is, let us suppose that the operator space
$E$ has Fourier type $p$ with respect to $G$. Then we wonder if
$\Lambda_{G,E,p}^1 (f) = \mathcal{F}_{G,E} (f)$ for all $f \in
L_E^p(G)$.

\begin{Lem} \label{LemmaType1}
$\|\widehat{f}\|_{\Le_E^{\infty}(\G)} \le \|f\|_{L_E^1(G)}$ for
all $f \in L_E^1(G)$.
\end{Lem}

\emph{Proof}. Since $E$ is an operator space we have $E \subset
\B(\Hi)$ for some Hilbert space $\Hi$. Hence, if $\mathbf{h} =
(h_1,h_2,\ldots,h_{d_{\pi}}) \in \ell_{\Hi}^2(d_{\pi})$, we can
write
\begin{eqnarray*}
\|\widehat{f}(\pi)\|_{S_{d_{\pi}}^{\infty}(E)} & = &
\|\widehat{f}(\pi)\|_{\B(\ell_{\Hi}^2(d_{\pi}))} \\ & \le &
\sup_{\|\mathbf{h}\|_{\ell_{\Hi}^2(d_{\pi})} \le 1} \Big(
\sum_{i=1}^{d_{\pi}} \big[ \int_G \Big\| f(g) \Big(
\sum_{j=1}^{d_{\pi}} \overline{\pi_{ji}(g)} h_j \Big) \Big\|_{\Hi}
d\mu(g) \big]^2 \Big)^{1/2}.
\end{eqnarray*}
Applying Minkowski inequality for integrals we get
$$\|\widehat{f}(\pi)\|_{S_{d_{\pi}}^{\infty}(E)} \le
\sup_{\|\mathbf{h}\|_{\ell_{\Hi}^2(d_{\pi})} \le 1} \int_G
\|f(g)\|_E \Big( \sum_{i=1}^{d_{\pi}} \big\| \sum_{j=1}^{d_{\pi}}
\overline{\pi_{ji}(g)} h_j \big\|_{\Hi}^2 \Big)^{1/2} d\mu(g).$$
Therefore we just need to check the inequality
$$\sup_{\|\mathbf{h}\|_{\ell_{\Hi}^2(d_{\pi})} \le 1} \Big(
\sum_{i=1}^{d_{\pi}} \big\| \sum_{j=1}^{d_{\pi}}
\overline{\pi_{ji}(g)} h_j \big\|_{\Hi}^2 \Big)^{1/2} \le 1$$ for
all $\pi \in \G$ and almost all $g \in G$. But this is a simple
consequence of the unitarity of $\pi(g)$ for any $g \in G$. This
completes the proof. \fin

\begin{Pro} Let $E$ be an operator space having Fourier type
$p$ with respect to $G$. Then $\Lambda_{G,E,p}^1 (f) =
\mathcal{F}_{G,E} (f)$ for all $f \in L_E^p(G)$.
\end{Pro}

\emph{Proof}. Let $\{f_n\}_{n=1}^{\infty} \subset L^p(G) \otimes
E$ be a sequence convergent to $f$ in the norm of $L_E^p(G)$.
Then, applying Lemma \ref{LemmaType1}, we have
\begin{eqnarray*} \lefteqn{\|\Lambda_{G,E,p}^1(f) -
\mathcal{F}_{G,E}(f)\|_{\Le_E^{\infty}(\G)}} \\ & \le &
\|\Lambda_{G,E,p}^1(f - f_n)\|_{\Le_E^{\infty}(\G)} +
\|\mathcal{F}_{G,E}(f_n - f)\|_{\Le_E^{\infty}(\G)} \\ & \le &
\|\Lambda_{G,E,p}^1(f - f_n)\|_{\Le_E^{p'}(\G)} +
\|\mathcal{F}_{G,E}(f_n - f)\|_{\Le_E^{\infty}(\G)} \\ & \le &
\Co_p^1(E,G) \|f - f_n\|_{L_E^p(G)} + \|f_n - f\|_{L_E^1(G)}
\\ & \le & \big( \Co_p^1(E,G) + 1 \big) \|f - f_n\|_{L_E^p(G)}
\end{eqnarray*}
The result follows by taking the limit in $n$. This completes the
proof. \fin

As it is well-known, every Banach space has Fourier type 1 in the
sense of Milman \cite{M}. In the following result, which extends
Lemma \ref{LemmaType1}, we show that every operator space has
Fourier type $1$.

\begin{Pro} \label{Type1} We have $\Co^1_1(E,G) = 1$ for every
pair $(E,G)$.
\end{Pro}

\emph{Proof}. Let us denote by $\max B$ the operator space which
results when we impose on the Banach space $B$ its max
quantization, see Chapter $3$ of \cite{ER2} for the details. Let
$E_1$ and $E_2$ be operator spaces. Then the natural
identification $cb((\max B) \otimes^{\wedge} E_1, E_2) \simeq
\mathcal{B}(B, cb(E_1, E_2))$, given by $\Psi(b \otimes e_1) =
\Phi(b)(e_1)$, is a completely isomorphic isomorphism. This
follows by the factorization $$cb((\max B) \otimes^{\wedge} E_1,
E_2) \simeq cb(\max B, cb(E_1, E_2)) \simeq \mathcal{B}(B, cb(E_1,
E_2))$$ which is composed of completely isometric isomorphisms,
see Chapters $3$ and $7$ of \cite{ER2}. Therefore, since the space
$L_E^1(G)$ can be rewritten as $\max L^1(G) \otimes^{\wedge} E$,
we get that $$\mathcal{C}_1^1(E,G) = \sup_{\|f\|_{L^1(G)} \le 1}
\|\widehat{f} \otimes \hspace{2pt} \cdot \hspace{2pt} \|_{cb(E,
\mathcal{L}_E^{\infty}(\G))} \le \sup_{\|f\|_{L^1(G)} \le 1}
\|\widehat{f}\|_{\mathcal{L}^{\infty}(\G)} = 1$$ by the
Hausdorff-Young inequality on compact groups, see \cite{Ku} or
Lemma \ref{H-Y} below. Recall that the supremum is attained taking
$f$ to be the constant function $1$. This completes the proof.
\fin

\begin{Rem}
\emph{There exists an alternative approach to this result using
similar arguments to those employed in the proof of Proposition
\ref{CotypeInfinite}.}
\end{Rem}

The following corollary exhibits the Fourier type as a stronger
condition on the pair $(E,G)$ as the exponent $p$ approaches $2$.
Its proof follows by means of Proposition \ref{Type1} and the
complex interpolation method.

\begin{Cor} Let $1 \le p_1 \le p_2 \le 2$ and assume that the
operator space $E$ has Fourier type $p_2$ with respect to $G$.
Then $E$ has Fourier type $p_1$ with respect to $G$. Moreover we
have $\Co_{p_1}^1(E,G) \le \Co_{p_2}^1(E,G)^{p_2'/p_1'}$.
\end{Cor}

A vector-valued version of the Riemann-Lebesgue lemma on compact
groups follows easily from Proposition \ref{Type1} and the scalar
result.

\begin{Cor} $\mathcal{F}_{G,E}
(L_E^1(G)) \subset C_0(\G,E)$ for every operator space $E$.
\end{Cor}

\section{Fourier cotype} \label{section5}

If $G$ is a locally compact abelian group, the Fourier inversion
theorem asserts that any $f \in L^1(G)$ such that $\widehat{f} \in
L^1(\G)$ can be recovered as $$f(g) = \,\
\widehat{\!\!\widehat{f}} (g^{-1}) \qquad \mbox{for almost every}
\ \ g \in G.$$ Furthermore, if $G$ is compact one can conclude
that the operators $\mathcal{F}_G^{-1}$ and $\mathcal{F}_{\G}$ are
essentially the same via the topological isomorphism from $G$ onto
its bidual, given by the Pontrjagin duality theorem. In the
vector-valued context this means that, in order to study the
operator $\mathcal{F}_{G,E}^{-1}$, it suffices to study the
Fourier transform $\mathcal{F}_{\G,E}$. For this reason we do not
find the concept of Fourier cotype in the commutative theory.
However, for a non-commutative compact group $G$, the Fourier
inversion theorem and the Pontrjagin duality theorem are no longer
valid since the dual object $\G$ is not even a group. These
consi\-derations explain why the study of the inverse operator
$\mathcal{F}_{G,E}^{-1}$ should not be a trivial consequence of
the analysis of the operator $\mathcal{F}_{G,E}$.

Let $E$ be an operator space and $1 \le p \le 2$. Arguing as in
section $\ref{section4}$, we can deduce the relation
$\mathcal{F}_{G,E}^{-1} (\Le^p(\G) \otimes E) \subset L^{p'}(G)
\otimes E$. This follows from Kunze's result for the inverse
Fourier transform on compact groups (see \cite{Ku} or Lemma
$\ref{H-Y}$ below). This motivates the following definition.

\begin{Def} \emph{Let $1 \le p \le 2$ and let $p'$ be the conjugate
exponent of $p$. We say that the operator space $E$ has
\textbf{Fourier cotype} $p'$ with respect to the compact group $G$
if the operator $\mathcal{F}_{G,E}^{-1}: \Le^p(\G) \otimes E
\rightarrow L^{p'}(G) \otimes E$ can be extended to a completely
bounded operator
$$\Lambda_{G,E,p'}^2: \Le_E^p(\G) \longrightarrow L_E^{p'}(G).$$
In that case, we shall denote by $\Co_{p'}^2(E,G)$ the $cb$ norm
of $\Lambda_{G,E,p'}^2$.}
\end{Def}

\begin{Rem}
\emph{Now it is obvious that, for compact abelian groups, the
notion of Fourier cotype is the completely bounded version of
Milman's notion of Fourier type with respect to the dual group
$\G$.}
\end{Rem}

Plancherel theorem for compact groups gives an explicit formula
for the action of $\mathcal{F}_G^{-1}$ on $\Le^2(\G)$ and, by the
natural embeddings, also on $\Le^p(\G)$ for $1 \le p \le 2$. It is
obvious that this formula remains valid if we take tensor
products. Namely, given $1 \le p \le 2$, the action of the
operator $\mathcal{F}_{G,E}^{-1}$ on $\Le^p(\G) \otimes E$ is
given by $$A \in \Le^p(\G) \otimes E \longmapsto \sum_{\pi \in \G}
d_{\pi} \mbox{tr} (A^{\pi} \pi(\cdot)) \in L^{p'}(G) \otimes E.$$
Therefore, if we want our definition of Fourier cotype to be
natural, we need affirmative answers for the following questions.

\begin{itemize}
\item[(a)] \emph{Does the operator $\Lambda_{G,E,p'}^2$ preserve
the given explicit formula?} That is, if the operator space $E$
has Fourier cotype $p'$ with respect to $G$, we ask whether for
all $A \in \Le_E^p(\G)$ we have $$\Lambda_{G,E,p'}^2(A) =
\sum_{\pi \in \G} d_{\pi} \mbox{tr} (A^{\pi} \pi(\cdot)).$$ If $A
\in \Le_E^p(\G)$, it has a countable support $I_A =
\{\pi_k\}_{k=1}^{\infty} \subset \G$. Then we define $A_n \in
\Le^p(\G) \otimes E$ by the relations $A_n^{\pi} = A^{\pi}$ if
$\pi = \pi_k$ for $1 \le k \le n$ and $A_n^{\pi} = 0$ otherwise.
Denoting by $$f = \sum_{\pi \in \G} d_{\pi} \mbox{tr} (A^{\pi}
\pi(\cdot)) \qquad \mbox{and} \qquad f_n = \sum_{\pi \in \G}
d_{\pi} \mbox{tr} (A_n^{\pi} \pi(\cdot))$$ we obtain that
\begin{eqnarray*}
\|\Lambda_{G,E,p'}^2(A) - f\|_{L_E^{p'}(G)} & \le &
\|\Lambda_{G,E,p'}^2(A - A_n)\|_{L_E^{p'}(G)} + \|f -
f_n\|_{L_E^{p'}(G)} \\ & \le & C_{p'}^2(E,G) \|A -
A_n\|_{\Le_E^p(\G)} + \|f - f_n\|_{L_E^{p'}(G)}.
\end{eqnarray*}
The first term of the sum is arbitrarily small as $n$ tends to
infinity. For the second term it is not difficult to check that
the sequence $\{f_n\}_{n=1}^{\infty}$ is Cauchy. Thus, replacing
this sequence if necessary by an appropriate subsequence, we can
assume $\|f_{n_2} - f_{n_1}\|_{L_E^{p'}(G)} < 2^{-m}$ for all
$n_1, n_2 \ge m$. Hence $$\|f - f_n\|_{L_E^{p'}(G)} \le
\sum_{k=n+1}^{\infty} \|f_k - f_{k-1}\|_{L_E^{p'}(G)} <
\sum_{k=n}^{\infty} \frac{1}{2^k}$$ and $\Lambda_{G,E,p'}^2(A) =
\displaystyle \sum_{\pi \in \G} d_{\pi} \mbox{tr} (A^{\pi}
\pi(\cdot))$ as we wanted.

\item[(b)] \emph{Does the operator $\Lambda_{G,E,p'}^2$ coincide
with the inverse of the vector-valued  Fourier transform?} That
is, if the operator space $E$ has Fourier cotype $p'$ with respect
to $G$, we ask whether for all $A \in \Le_E^p(\G)$ we have
$$\mathcal{F}_{G,E} \circ \Lambda_{G,E,p'}^2 (A) = A.$$ Using the
same notation as above we just need to see that $\widehat{f} = A$.
Given $\pi \in \G$ we take $n_{\pi}$ to be the smallest positive
integer satisfying $\pi \neq \pi_k$ for $k \ge n_{\pi}$. Then it
is obvious that $\widehat{f}(\pi) - A^{\pi} = (\widehat{f} -
\widehat{f}_n) (\pi)$ for all $n \ge n_{\pi}$ and therefore it is
enough to estimate the entries of that matrix. Namely,
\begin{eqnarray*}
\big\| \big( (\widehat{f} - \widehat{f}_n) (\pi) \big)_{ij}
\big\|_E & \le & \int_G \|(f -f_n)(g)\|_E |\pi_{ji}(g)| \, d\mu(g)
\\ & \le & \|f - f_n\|_{L_E^{p'}(G)} \, \|\pi_{ji}\|_{L^p(G)} \\
& \le & \Co_{p'}^2 (E,G) \, \|A - A_n\|_{\Le_E^p(\G)}
\end{eqnarray*}
which is arbitrarily small for large $n$.
\end{itemize}

\begin{Pro} \label{CotypeInfinite} We have $\Co^2_{\infty}(E,G) =
1$ for every pair $(E,G)$.
\end{Pro}

\emph{Proof}. By property $1$ of Theorem \ref{Schatten} and a
density argument, we have to see that for all $n \ge 1$, any
family of vectors $A_{ij} \in \Le^1(\G) \otimes E$ and almost all
$g \in G$, we have $$\Big\| \Big( \,\ \sum_{\pi \in \G} d_{\pi}
\mbox{tr} (A_{ij}^{\pi} \pi(g)) \,\ \Big) \Big\|_{S_n^{\infty}(E)}
\le \Big\| \Big( \,\ A_{ij} \,\ \Big)
\Big\|_{S_n^{\infty}(\Le_E^1(\G))}.$$ If we consider a vector $A
\in \Le_E^1(\G)$ as an element of $cb(\Le^{\infty}(\G), E)$ by the
relation $$B \in \Le^{\infty}(\G) \longmapsto \sum_{\pi \in \G}
d_{\pi} \mbox{tr} (A^{\pi} B^{\pi}) \in E$$ then it is easy to see
that, for $B_g \in \Le^{\infty}(\G)$ defined by $B_g^{\pi} =
\pi(g)$, we have
\begin{eqnarray*}
\Big\| \Big( \,\ \sum_{\pi \in \G} d_{\pi} \mbox{tr} (A_{ij}^{\pi}
\pi(g)) \,\ \Big) \Big\|_{S_n^{\infty}(E)} & = & \Big\| \Big( \,\
\sum_{\pi \in \G} d_{\pi} \mbox{tr} (A_{ij}^{\pi} \cdot) (B_g) \,\
\Big) \Big\|_{S_n^{\infty}(E)} \\ & \le & \Big\| \Big( \,\
\sum_{\pi \in \G} d_{\pi} \mbox{tr} (A_{ij}^{\pi} \cdot) \,\ \Big)
\Big\|_{cb(\Le^{\infty}(\G),S_n^{\infty}(E))} \\ & \le & \Big\|
\Big( \,\ A_{ij} \,\ \Big) \Big\|_{S_n^{\infty}(\Le^1(\G)
\otimes^{\wedge} E)}
\end{eqnarray*}
where the last inequality follows from the complete contraction
given by $\Le^1(\G) \otimes^{\wedge} E \rightarrow \Le^1(\G)
\otimes_{\mbox{{\tiny min}}} E \rightarrow
cb(\Le^{\infty}(\G),E)$. Finally, we get the desired relation by
Corollary \ref{L-one}. We have shown that $\Co_{\infty}^2(E,G) \le
1$. The reverse inequality follows from Corollary \ref{LowerBound}
below. \fin

\begin{Cor} Let $1 \le p_1 \le p_2 \le 2$ and assume that the
operator space $E$ has Fourier cotype $p'_2$ with respect to $G$.
Then $E$ has Fourier cotype $p'_1$ with respect to $G$. Moreover
we have $\Co_{p'_1}^2(E,G) \le \Co_{p'_2}^2(E,G)^{p_2'/p_1'}$.
\end{Cor}

\section{Duality, $cb$ distance and some other topics}
\label{section6}

\setcounter{The}{0}

Let $E$ be an operator space. The aim of this section is to study
the Fourier type and cotype of some operator spaces related to
$E$. We begin by stating the scalar-valued Hausdorff-Young
inequality. Recall that we write $\,\ \bcheck{A}$ to denote
$\mathcal{F}_G^{-1}(A)$.

\begin{Lem} [Hausdorff-Young inequality] \label{H-Y} Let $1 \le p
\le 2$ and let $p'$ be the conjugate exponent of $p$:
\begin{itemize}
\item[$1.$ ] If $f \in L^p(G)$, then $\widehat{f} \in \Le^{p'}(\G)$
and $\|\mathcal{F}_G\|_{cb(L^p(G),\Le^{p'}(\G))} =1$.
\item[$2.$ ] If $A \in \Le^p(\G)$, then $\,\ \bcheck{A} \in L^{p'}(G)$
and $\|\mathcal{F}_G^{-1}\|_{cb(\Le^p(\G),L^{p'}(G))} =1$.
\end{itemize}
\end{Lem}

Note that this statement of the inequality goes a bit further than
Kunze's original result since we are asserting that the Fourier
transform is not only bounded but completely bounded. The proof is
straightforward, first one checks that $\mathcal{F}_G$ is a
complete contraction from $L^1(G)$ into $\Le^{\infty}(\G)$. But,
since $L^1(G)$ is equipped with its max operator space structure,
the $cb$ norm coincides with the operator norm, see Chapter $3$
\cite{ER2}. The same argument works to see that the inverse
Fourier transform is a complete contraction from $\Le^1(\G)$ to
$L^{\infty}(G)$, now $L^{\infty}(G)$ is equipped with its min
operator space structure. These facts can also be justified as
simple consequences of Propositions \ref{Type1} and
\ref{CotypeInfinite}. Second, from the Plancherel theorem for
compact groups, it is easy to check that $\mathcal{F}_G$ is a
complete isometric isomorphism from $L^2(G)$ onto $\Le^2(\G)$. By
complex interpolation the general case is obtained and, the fact
that the $cb$ norm of the Fourier transform is not smaller than
$1$ for any $1 \le p \le 2$ can be checked by testing with the
constant function $1$.

\subsection{Basic results}

We begin by the simplest case. Namely, the Fourier type and cotype
of the subspaces of $E$. The following result is a trivial
consequence of property $1$ of Theorem \ref{Schatten}.

\begin{Pro} \label{Subspace} Let $F$ be a closed subspace of $E$,
then we have the estimates $\Co_p^1(F,G) \le \Co_p^1(E,G)$ and
$\Co_{q'}^2(F,G) \le \Co_{q'}^2(E,G)$ for any $1 \le p, q \le 2$.
\end{Pro}

\begin{Cor} \label{LowerBound} $\Co_p^1(E,G) \ge 1$ and
$\Co_{q'}^2(E,G) \ge 1$ for any $1 \le p, q \le 2$.
\end{Cor}

Now we consider complex interpolation of operator spaces. The
proof of the next result is also straightforward.

\begin{Pro} Let $1 \le p_0, p_1 \le 2$ and assume that
$\{E_0, E_1\}$ is compatible for complex interpolation. Then
$\Co_{p_{\theta}}^1(E_{\theta},G) \le \Co_{p_0}^1(E_0,G)^{1 -
\theta} \Co_{p_1}^1(E_1,G)^{\theta}$ for $p_{\theta}^{-1} = (1 -
\theta) p_0^{-1} + \theta p_1^{-1}$. A similar result holds for
the Fourier cotype.
\end{Pro}

\subsection{Duality}

The following theorem can be rephrased by saying that Fourier type
and cotype are dual notions.

\begin{The} \label{Duality} Let $E$ be an operator space, $1
\le p \le 2$ and $p'$ its conjugate exponent. Then
\begin{itemize}
\item[$1.$ ] $E$ has Fourier type $p$ with respect to a
compact group $G$ if and only if $E^{\star}$ has Fourier cotype
$p'$ with respect to $G$.
\item[$2.$ ] $E$ has Fourier cotype $p'$ with respect to a
compact group $G$ if and only if $E^{\star}$ has Fourier type $p$
with respect to $G$.
\end{itemize}
Moreover, we have $\Co_p^1(E,G) = \Co_{p'}^2(E^{\star},G)$ and
$\Co_p^1(E^{\star},G) = \Co_{p'}^2(E,G)$.
\end{The}

\emph{Proof}. We just prove the equality $\Co_p^1(E,G) =
\Co_{p'}^2(E^{\star},G)$ since the proof of the second identity is
essentially the same. The case $p = 1$ follows from Propositions
\ref{Type1} and \ref{CotypeInfinite}, thus we assume that $1 < p
\le 2$.

\

\emph{Step} 1. $\Co_p^1(E,G) \ge \Co_{p'}^2(E^{\star},G)$. By a
density argument and property $1$ of Theorem \ref{Schatten} we
just need to check that the following inequality holds $$\Big\|
\Big( \,\ \sum_{\pi \in \G} d_{\pi} \mbox{tr} (A_{ij}^{\pi}
\pi(\cdot)) \,\ \Big) \Big\|_{S_n^{p'}(L_{E^{\star}}^{p'}(G))} \le
\Co_p^1(E,G) \ \ \Big\| \Big( \,\ A_{ij} \,\ \Big)
\Big\|_{S_n^{p'}(\Le_{E^{\star}}^p(\G))}$$ for any family $A_{ij}
\in \Le^p(\G) \otimes E^{\star}$ ($1 \le i,j \le n$) and all $n
\ge 1$. But we have the completely isometric isomorphism
$S_n^{p'}(L_{E^{\star}}^{p'}(G)) \simeq
L_{S_n^p(E)^{\star}}^{p'}(G)$. So for all $\varepsilon
> 0$ there exists $f^{\varepsilon} \in L_{S_n^p(E)}^p(G)$ of norm $1$
such that
\begin{eqnarray*}
\lefteqn{\Big\| \Big( \,\ \sum_{\pi \in \G} d_{\pi} \mbox{tr}
(A_{ij}^{\pi} \pi(\cdot)) \,\ \Big)
\Big\|_{S_n^{p'}(L_{E^{\star}}^{p'}(G))}} \\ & \le & (1 +
\varepsilon) \, \Big| \int_G \mbox{tr} \Big[ \Big( \,\ \sum_{\pi
\in \G} d_{\pi} \mbox{tr} (A_{ij}^{\pi} \pi(g)) \,\ \Big) \Big(
\,\ f_{ij}^{\varepsilon}(g) \,\ \Big) \Big] d \mu(g) \Big|
\end{eqnarray*}
and where $f_{ij}^{\varepsilon}$, the entries of
$f^{\varepsilon}$, belong to $L^p(G) \otimes E$. If we denote by
$\mathcal{I}$ the integral over $G$ written above, then we would
like to prove that $$\mathcal{I} = \sum_{i,j=1}^n \sum_{\pi \in
\G} d_{\pi} \ \ \int_G \big\langle \mbox{tr} (A_{ij}^{\pi}
\pi(g)), f_{ji}^{\varepsilon}(g) \big\rangle d \mu(g).$$ Taking
into account that $A_{ij} \in \Le^p(\G) \otimes E^{\star}$ and
$f_{ij}^{\varepsilon} \in L^p(G) \otimes E$ it suffices to show
that the expressions
\begin{eqnarray*}
\mathcal{I}_1 & = & \int_G \sum_{\pi \in \G} d_{\pi} \mbox{tr}
(A^{\pi} \pi(g)) \,\ f(g) \,\ d \mu(g) \\ \mathcal{I}_2 & = &
\sum_{\pi \in \G} d_{\pi} \int_G \mbox{tr} (A^{\pi} \pi(g)) \,\
f(g) \,\ d \mu(g)
\end{eqnarray*}
coincide for all $A \in \Le^p(\G)$ and all $f \in L^p(G)$. But
this is an easy computation that we leave to the reader. In
summary we obtain
\begin{eqnarray*}
\lefteqn{\Big\| \Big( \,\ \sum_{\pi \in \G} d_{\pi} \mbox{tr}
(A_{ij}^{\pi} \pi(\cdot)) \,\ \Big)
\Big\|_{S_n^{p'}(L_{E^{\star}}^{p'}(G))}} \\ & \le & (1 +
\varepsilon) \, \Big| \sum_{i,j=1}^n \sum_{\pi \in \G} d_{\pi}
\int_G \big\langle \mbox{tr} (A_{ij}^{\pi} \pi(g)),
f_{ji}^{\varepsilon}(g) \big\rangle d \mu(g) \Big| \\ & = & (1 +
\varepsilon) \, \Big| \sum_{i,j=1}^n \sum_{\pi \in \G} d_{\pi}
\mbox{tr}\big( \langle A_{ij}^{\pi}, \widehat{\tau
(f_{ji}^{\varepsilon})}(\pi) \rangle \big) \Big| \\ & = & (1 +
\varepsilon) \, \Big| \mbox{tr} \Big[ \Big( \,\ A_{ij} \,\ \Big)
\Big( \,\ \widehat{\tau (f_{ij}^{\varepsilon})} \,\ \Big) \Big]
\Big|
\end{eqnarray*}
where $\tau(f)(g) = f(g^{-1})$. This step is concluded by the
following inequality
\begin{eqnarray*}
\Big| \mbox{tr} \Big[ \Big( \,\ A_{ij} \,\ \Big) \Big( \,\
\widehat{\tau (f_{ij}^{\varepsilon})} \,\ \Big) \Big] \Big| & \le
& \Big\| \Big( \,\ A_{ij} \,\ \Big)
\Big\|_{S_n^{p'}(\Le_{E^{\star}}^p(\G))} \ \ \Big\| \Big( \,\
\widehat{\tau (f_{ij}^{\varepsilon})} \,\ \Big)
\Big\|_{S_n^p(\Le_E^{p'}(\G))} \\ & \le & \Co_p^1(E,G) \ \ \Big\|
\Big( \,\ A_{ij} \,\ \Big)
\Big\|_{S_n^{p'}(\Le_{E^{\star}}^p(\G))}
\|f^{\varepsilon}\|_{S_n^p(L_E^p(G))}.
\end{eqnarray*}

\emph{Step} 2. $\Co_p^1(E,G) \le \Co_{p'}^2(E^{\star},G)$. By the
same reasons given in Step $1$, it suffices to check that $$\Big\|
\Big( \,\ \widehat{f}_{ij} \,\ \Big)
\Big\|_{S_n^{p'}(\Le_E^{p'}(\G))} \le \Co_{p'}^2(E^{\star},G) \ \
\Big\| \Big( \,\ f_{ij} \,\ \Big) \Big\|_{S_n^{p'}(L_E^p(G))}$$
for any family $f_{ij} \in L^p(G) \otimes E$ ($1 \le i,j \le n$)
and all $n \ge 1$. Given $\varepsilon > 0$, the complete isometry
$$S_n^{p'}(\Le_E^{p'}(\G)) \simeq \Le_{S_n^{p'}(E)}^{p'}(\G)$$
provides the existence of $\mathbf{A}^{\varepsilon} \in
\Le_{S_n^p(E^{\star})}^p(\G)$ of norm $1$ such that $$\Big\| \Big(
\,\ \widehat{f}_{ij} \,\ \Big) \Big\|_{S_n^{p'}(\Le_E^{p'}(\G))}
\le (1 + \varepsilon) \, \Big| \sum_{\pi \in \G} d_{\pi} \mbox{tr}
\Big[ \Big( \,\ A_{ij}^{\varepsilon, \pi} \,\ \Big) \Big( \,\
\widehat{f}_{ij}(\pi) \,\ \Big) \Big] \Big|$$ and where
$A_{ij}^{\varepsilon}$, the entries of $\mathbf{A}^{\varepsilon}$,
belong to $\Le^p(\G) \otimes E^{\star}$. If $\mathcal{S}$ denotes
the sum written above, then we can argue as in Step $1$ to obtain
$$\mathcal{S} = \sum_{i,j=1}^n \int_G \big\langle \sum_{\pi \in
\G} d_{\pi} \mbox{tr} (A_{ij}^{\varepsilon, \pi} \pi(g)^{\star}),
f_{ji}(g) \big\rangle d \mu(g).$$ Therefore
\begin{eqnarray*}
\Big\| \Big( \,\ \widehat{f}_{ij} \,\ \Big)
\Big\|_{S_n^{p'}(\Le_E^{p'}(\G))} & \le & (1 + \varepsilon) \,
\Big| \sum_{i,j=1}^n \int_G \big\langle
\mathcal{F}_{G,E^{\star}}^{-1} (A_{ij}^{\varepsilon})(g),
f_{ji}(g^{-1}) \big\rangle d \mu(g) \Big| \\ & = & (1 +
\varepsilon) \, \Big| \mbox{tr} \Big[ \Big( \,\
\mathcal{F}_{G,E^{\star}}^{-1} (A_{ij}^{\varepsilon}) \,\ \Big)
\Big( \,\ \tau (f_{ij}) \,\ \Big) \Big] \Big| \\ & \le & (1 +
\varepsilon) \,\ \big\| \mathcal{F}_{G,E^{\star}}^{-1}
(\mathbf{A}^{\varepsilon}) \big\|_{S_n^p(L_E^p(G)^{\star})} \Big\|
\Big( \,\ f_{ij} \,\ \Big) \Big\|_{S_n^{p'}(L_E^p(G))} \\ & = & (1
+ \varepsilon) \,\ \big\| \mathcal{F}_{G,E^{\star}}^{-1}
(\mathbf{A}^{\varepsilon}) \big\|_{S_n^p(L_{E^{\star}}^{p'}(G))}
\Big\| \Big( \,\ f_{ij} \,\ \Big) \Big\|_{S_n^{p'}(L_E^p(G))} \\ &
\le & (1 + \varepsilon) \,\ \Co_{p'}^2(E^{\star},G) \,\ \Big\|
\Big( \,\ f_{ij} \,\ \Big) \Big\|_{S_n^{p'}(L_E^p(G))}
\end{eqnarray*}
The proof is completed by taking $\varepsilon$ arbitrarily small.
\fin

\begin{Rem}
\emph{There exists another possible approach to this result.
Namely, $E$ has Fourier type $p$ if and only if the corresponding
Fourier transform ope\-rator is completely bounded. But then, by
Proposition $3.2.2$ of \cite{ER2}, the adjoint operator is also
completely bounded with the same $cb$ norm. Moreover, it can be
checked that the adjoint coincides with the inverse of the Fourier
transform for functions taking values in $E^{\star}$. This gives
the first equality of Theorem \ref{Duality}. The second equality
follows in a similar fashion.}
\end{Rem}

\begin{Cor} $\Co_p^1(E,G) = \Co_p^1(E^{\star
\star},G)$ and $\Co_{p'}^2(E,G) = \Co_{p'}^2(E^{\star \star},G)$.
\end{Cor}

\subsection{The $cb$ distance}

There exists a natural analog in the category of operator spaces
of the Banach-Mazur distance, due to Pisier. It is called the $cb$
distance and it is defined by $$d_{cb} (E_1,E_2) = \inf
\{\|u\|_{cb(E_1,E_2)} \|u^{-1}\|_{cb(E_2,E_1)}\}$$ where the
infimum runs over all complete isomorphisms $u: E_1 \rightarrow
E_2$.

\begin{The} \label{Local} Let $E_1$ and $E_2$ be operator spaces
and let $G$ be a compact group. Then the following inequalities
hold for $1 \le p \le 2$
\begin{eqnarray*}
\Co_p^1(E_2,G) & \le & d_{cb}(E_1,E_2) \,\ \Co_p^1(E_1,G) \\
\Co_{p'}^2(E_2,G) & \le & d_{cb}(E_1,E_2) \,\ \Co_{p'}^2(E_1,G) \\
\Co_p^1(E_2,G) & \le & d_{cb}(E_1,E_2^{\star}) \,\
\Co_{p'}^2(E_1,G) \\ \Co_{p'}^2(E_2,G) & \le &
d_{cb}(E_1,E_2^{\star}) \,\ \Co_p^1(E_1,G).
\end{eqnarray*}
\end{The}

\emph{Proof}. The last two inequalities follow from the first two
ones plus duality. Let us assume that the first inequality holds,
then the second inequality is also an immediate consequence of
Theorem \ref{Duality}
\begin{eqnarray*}
\Co_{p'}^2(E_2,G) & = & \Co_p^1(E_2^{\star},G) \,\ \le \,\
d_{cb}(E_1^{\star},E_2^{\star}) \,\ \Co_p^1(E_1^{\star},G) \\ & =
& d_{cb}(E_1^{\star},E_2^{\star}) \,\ \Co_{p'}^2(E_1,G) \,\ = \,\
d_{cb}(E_1,E_2) \,\ \Co_{p'}^2(E_1,G)
\end{eqnarray*}
where we have applied the identity $\|u\|_{cb(E_1,E_2)} =
\|u^{\star}\|_{cb(E_2^{\star},E_1^{\star})}$ to justify the
equality $d_{cb}(E_1^{\star},E_2^{\star}) = d_{cb}(E_1,E_2)$, see
\cite{ER2} for details. Therefore we will be done if we prove the
validity of the first inequality. For that it suffices to see
$$\Big\| \Big( \,\ \widehat{f}_{ij} \,\ \Big)
\Big\|_{S_n^{p'}(\Le_{E_2}^{p'}(\G))} \le \|u\|_{cb} \,\
\|u^{-1}\|_{cb} \,\ \Co_p^1(E_1,G) \,\ \Big\| \Big( \,\ f_{ij} \,\
\Big) \Big\|_{S_n^{p'}(L_{E_2}^p(G))}$$ for any family $f_{ij} \in
L^p(G) \otimes E_2$ ($1 \le i,j \le n$), any complete isomorphism
$u: E_1 \rightarrow E_2$ and all $n \ge 1$. But
\begin{eqnarray*}
\lefteqn{\Big\| \Big( \,\ \widehat{f}_{ij} \,\ \Big)
\Big\|_{S_n^{p'}(\Le_{E_2}^{p'}(\G))} = \Big( \sum_{\pi \in \G}
d_{\pi} \Big\| \Big( \,\ \widehat{f}_{ij}(\pi) \,\ \Big)
\Big\|_{S_{d_{\pi}n}^{p'}(E_2)}^{p'} \Big)^{1/p'}} \\ & \le &
\|u\|_{cb} \,\ \Big( \sum_{\pi \in \G} d_{\pi} \Big\|
(I_{M_{d_{\pi}n}} \otimes u^{-1}) \Big( \,\ \widehat{f}_{ij}(\pi)
\,\ \Big) \Big\|_{S_{d_{\pi}n}^{p'}(E_1)}^{p'} \Big)^{1/p'} \\ & =
& \|u\|_{cb} \,\ \Big\| \Big( \,\ \mathcal{F}_{G, E_1} \big(
(I_{L^p(G)} \otimes u^{-1}) (f_{ij}) \big) \,\ \Big)
\Big\|_{\Le_{S_n^{p'}(E_1)}^{p'}(\G)}
\\ & \le & \|u\|_{cb} \,\ \Co_p^1(E_1,G) \,\ \Big\| \Big( \,\
(I_{L^p(G)} \otimes u^{-1}) (f_{ij}) \,\ \Big)
\Big\|_{S_n^{p'}(L_{E_1}^p(G))} \\ & \le & \|u\|_{cb} \,\
\Co_p^1(E_1,G) \,\ \|I_{L^p(G)} \otimes u^{-1}\|_{cb} \,\ \Big\|
\Big( \,\ f_{ij} \,\ \Big) \Big\|_{S_n^{p'}(L_{E_2}^p(G))} \\ & =
& \|u\|_{cb} \,\ \|u^{-1}\|_{cb} \,\ \Co_p^1(E_1,G) \,\ \Big\|
\Big( \,\ f_{ij} \,\ \Big) \Big\|_{S_n^{p'}(L_{E_2}^p(G))}.
\end{eqnarray*}
This completes the proof. \fin

We recall here that, if $\mathcal{OS}_n$ denote the class of all
$n$-dimensional ope\-rator spaces, Pisier proved the estimate
$d_{cb}(E,OH_n) \le \sqrt{n}$ for any operator space $E \in
\mathcal{OS}_n$. Here $OH_n$ denotes the $n$-dimensional operator
Hilbert space $OH$, see \cite{P1}. Therefore, by taking $E_1 =
l^2(n)$ in Theorem \ref{Local} and invoking the results of the
next section, we get the following result.

\begin{Cor}
We have $\Co_2^1(E,G), \,\ \Co_2^2(E,G) \le \sqrt{n}$ for any $E
\in \mathcal{OS}_n$.
\end{Cor}

\section{Basic examples}
\label{section7}

We study here the Fourier type and cotype of Lebesgue spaces,
Schatten classes and their vector-valued versions. We start with
the statement of some inequalities of Minkowski type in the
operator space setting.

If $1 \le p_1 \le p_2 \le \infty$ and the measure spaces
$(\Omega_1, \mathcal{M}_1, \nu_1)$, $(\Omega_2, \mathcal{M}_2,
\nu_2)$ are $\sigma$-finite, then the classical Minkowski
inequality for integrals asserts that the natural map
$$L_{L^{p_2}(\Omega_2)}^{p_1}(\Omega_1) \longrightarrow
L_{L^{p_1}(\Omega_1)}^{p_2}(\Omega_2)$$ is contractive. The same
happens if our functions $f: \Omega_1 \times \Omega_2 \rightarrow
E$ take values in a Banach space $E$. We are interested in the
complete boundedness of this operator and some others in which we
replace the Lebesgue spaces $L^p(\Omega)$ by the Schatten classes
$S_n^p$. For this purpose, by complex interpolation, it suffices
to check the cases $p_1 = p_2$ and $(p_1,p_2) = (1,\infty)$. The
first case follows from the Fubini type results stated in Theorem
\ref{Schatten}. The second case reduces to see that the natural
map $E_1 \otimes^{\wedge} (E_2 \otimes_{\mbox{{\tiny min}}} E_3)
\rightarrow (E_1 \otimes^{\wedge} E_2) \otimes_{\mbox{{\tiny
min}}} E_3$ is a complete contraction. The proof of this result
can be found in Theorem $8.1.10$ of \cite{ER2}. In summary we can
state the following results.

\begin{The} [Quantized Minkowski inequalities] \label{Minkowski}
Let us consider an operator space $E$ and let $1 \le p_1 \le p_2
\le \infty$.
\begin{itemize}
\item[$1.$ ] \emph{Lebesgue spaces}. Let $(\Omega_1,
\mathcal{M}_1, \nu_1)$ and $(\Omega_2, \mathcal{M}_2, \nu_2)$ be
$\sigma$-finite measure spaces. Then the following natural map is
a complete contraction $$L_{L_E^{p_2}(\Omega_2)}^{p_1}(\Omega_1)
\longrightarrow L_{L_E^{p_1}(\Omega_1)}^{p_2}(\Omega_2).$$

\item[$2.$ ] \emph{Schatten classes}. Let $k_1,
k_2 \ge1$, then the following natural map is a complete
contraction $$S_{k_1}^{p_1}(S_{k_2}^{p_2}(E)) \longrightarrow
S_{k_2}^{p_2}(S_{k_1}^{p_1}(E)).$$

\item[$3.$ ] \emph{Combined results}. Let $(\Omega, \mathcal{M},
\nu)$ be a measure space and $k \ge 1$. Then the following natural
maps are complete contractions $$S_k^{p_1}(L_E^{p_2}(\Omega))
\longrightarrow L_{S_k^{p_1}(E)}^{p_2}(\Omega) \ \ \ \ \
\mbox{and} \ \ \ \ \ L_{S_k^{p_2}(E)}^{p_1}(\Omega)
\longrightarrow S_k^{p_2}(L_E^{p_1}(\Omega)).$$
\end{itemize}
\end{The}

\begin{Rem}
\emph{The arguments sketched above in order to prove Theorem
\ref{Minkowski} need extra hypotheses. A different proof, without
those unnecessary hypotheses, can be found in the Thesis \cite{Pa}
of the second-named author.}
\end{Rem}

In the study of the Fourier type of a Banach space with respect to
a locally compact abelian group, Andersson \cite{A} gave the
following version of Minkowski inequality for regular measures. We
recall that our notion of regular measure is the same as the one
given in \cite{C}.

\begin{Pro} [Andersson] \label{Andersson} Let $1 \le p_1 \le p_2
< \infty$ and assume that $(\Omega_1, \mathcal{M}_1, \nu_1)$ and
$(\Omega_2, \mathcal{M}_2, \nu_2)$ are regular measure spaces. Let
us denote by $H$ the space of functions $f: \Omega_1 \times
\Omega_2 \rightarrow \C$ such that $|f|$ is bounded lower
semicontinuous and $\|f_{\omega_2}\|_{L^{p_1}(\Omega_1)}$ is
bounded in $\Omega_2$. Then the following natural map is
contractive $$L_{L^{p_2}(\Omega_2)}^{p_1}(\Omega_1) \cap H
\longrightarrow L_{L^{p_1}(\Omega_1)}^{p_2}(\Omega_2) \cap H.$$
\end{Pro}

Let us note that, if $(\Omega, \mathcal{M}, \nu)$ denotes a
regular measure space and we take $\Omega_1 = G$ and $\Omega_2 =
\Omega$ in Proposition \ref{Andersson}, then the space $C_c(G
\times \Omega)$ of continuous functions with compact support,
defined on $G \times \Omega$ and with values in $\C$, is contained
in the space $H$. Hence, by the density of $C_c(G \times \Omega)$
in $L_{L^{p_2}(\Omega)}^{p_1}(G)$ and
$L_{L^{p_1}(G)}^{p_2}(\Omega)$, we deduce that the natural map
$$L_{L^{p_2}(\Omega)}^{p_1}(G) \longrightarrow
L_{L^{p_1}(G)}^{p_2}(\Omega)$$ is a contraction whenever $1 \le
p_1 \le p_2 < \infty$. Then, by the same arguments that we gave in
the proof of Theorem \ref{Minkowski}, we conclude that we have in
fact a complete contraction. Furthermore, the same happens if we
take $\Omega_1 = \Omega$ and $\Omega_2 = G$. Therefore we have
shown the validity of the following result, which we enunciate for
vector-valued functions since its proof is analogous.

\begin{Lem} \label{Regular} Let $1 \le p_1 \le p_2
< \infty$ and assume that $E$ is an operator space, $G$ is a
compact group and $(\Omega, \mathcal{M}, \nu)$ is a regular
measure space. Then the following natural maps are complete
contractions $$L_{L_E^{p_2}(\Omega)}^{p_1}(G) \longrightarrow
L_{L_E^{p_1}(G)}^{p_2}(\Omega) \qquad \mbox{and} \qquad
L_{L_E^{p_2}(G)}^{p_1}(\Omega) \longrightarrow
L_{L_E^{p_1}(\Omega)}^{p_2}(G).$$
\end{Lem}

\begin{The} \label{BochnerLebesgue} Let $1 \le p, q \le 2$ and
assume that $E$ is an operator space having Fourier type $p$ and
Fourier cotype $q'$ with respect to a compact group $G$. Let
$(\Omega, \mathcal{M}, \nu)$ be a regular or $\sigma$-finite
measure space. Then
\begin{itemize}
\item[$1.$ ] $L_E^r(\Omega)$ has Fourier type $p$ with respect
to $G$ for all $p \le r \le p'$.
\item[$2.$ ] $L_E^s(\Omega)$ has Fourier cotype $q'$ with
respect to $G$ for all $q \le s \le q'$.
\end{itemize}
Moreover, $\Co_p^1(L_E^r(\Omega),G) = \Co_p^1(E,G)$ and
$\Co_{q'}^2(L_E^s(\Omega),G) = \Co_{q'}^2(E,G)$.
\end{The}

\emph{Proof}. We start by proving the relation
$\Co_p^1(L_E^r(\Omega),G) = \Co_p^1(E,G)$. For $p = 1$ we just
need to apply Proposition \ref{Type1}. Thus we assume that $1 < p
\le 2$. The inequality $\Co_p^1(L_E^r(\Omega),G) \ge \Co_p^1(E,G)$
follows from Proposition \ref{Subspace} and, by complex
interpolation, it then suffices to see that
$\Co_p^1(L_E^r(\Omega),G) \le \Co_p^1(E,G)$ for $r = p$ and $r =
p'$. For $r = p$ we observe that the natural map
$$L_{\Le_E^{p'}(\G)}^p(\Omega) \longrightarrow
\Le_{L_E^p(\Omega)}^{p'}(\G)$$ is a complete contraction. The
proof of this fact is similar to that of Theorem \ref{Minkowski}.
In particular we have
\begin{eqnarray*}
\Big\| \Big( \,\ \widehat{f}_{ij} \,\ \Big)
\Big\|_{S_n^p(\Le_{L_E^p(\Omega)}^{p'}(\G))} & \le & \Big\| \Big(
\,\ \widehat{f}_{ij} \,\ \Big)
\Big\|_{L_{S_n^p(\Le_E^{p'}(\G))}^p(\Omega)} \\ & \le &
\Co_p^1(E,G) \,\ \Big\| \Big( \,\ f_{ij} \,\ \Big)
\Big\|_{S_n^p(L_{L_E^p(\Omega)}^p(G))}.
\end{eqnarray*}
For $r = p'$ we use Theorem \ref{Minkowski} or Lemma
\ref{Regular}, depending on the measure space $(\Omega,
\mathcal{M}, \nu)$, to get the desired relation
\begin{eqnarray*}
\Big\| \Big( \,\ \widehat{f}_{ij} \,\ \Big)
\Big\|_{S_n^{p'}(\Le_{L_E^{p'}(\Omega)}^{p'}(\G))} & \le &
\Co_p^1(E,G) \,\ \Big\| \Big( \,\ f_{ij} \,\ \Big)
\Big\|_{L_{S_n^{p'}(L_E^p(G))}^{p'}(\Omega)} \\ & \le &
\Co_p^1(E,G) \,\ \Big\| \Big( \,\ f_{ij} \,\ \Big)
\Big\|_{S_n^{p'}(L_{L_E^{p'}(\Omega)}^p(G))}.
\end{eqnarray*}
The proof of the inequality $\Co_{q'}^2(L_E^s(\Omega),G) \le
\Co_{q'}^2(E,G)$ is analogous. \fin

\begin{Rem}
\emph{The proof of Theorem \ref{BochnerLebesgue} for scalar-valued
Lebesgue spaces is much simpler. Namely, one only has to see that
$L^2(\Omega)$ has Fourier type $2$ and then the result follows by
duality and complex interpolation with the trivial cases $p = 1$
and $p = \infty$. But the case $p = 2$ is a simple consequence of
Plancherel theorem on compact groups.}
\end{Rem}

It is well known that the dual of $L_E^p(\Omega)$ is not in
general $L_{E^{\star}}^{p'}(\Omega)$. However it is so when the
dual $E^{\star}$ possesses the Radon Nikodym property $RNP$. In
\cite{P2} Pisier developed an operator space version of the Radon
Nikodym property which he called $ORNP$. The following corollary,
which is a very simple consequence of Theorems \ref{Duality} and
\ref{BochnerLebesgue}, shows that both spaces have the same
Fourier type and cotype even if $E^{\star}$ does not satisfy the
$ORNP$.

\begin{Cor} Let $1 \le p, q \le 2$ and
assume that $E$ is an operator space having Fourier type $p$ and
Fourier cotype $q'$ with respect to a compact group $G$. Let
$(\Omega, \mathcal{M}, \nu)$ be a regular or $\sigma$-finite
measure space. Then we have
\begin{itemize}
\item[$1.$ ] $\Co_q^1(L_E^s(\Omega)^{\star},G) =
\Co_q^1(L_{E^{\star}}^{s'}(\Omega),G)$ for all $q \le s \le q'$.
\item[$2.$ ] $\Co_{p'}^2(L_E^r(\Omega)^{\star},G) =
\Co_{p'}^2(L_{E^{\star}}^{r'}(\Omega),G)$ for all $p \le r \le
p'$.
\end{itemize}
\end{Cor}

We now study the Fourier type and cotype of Schatten classes. We
will denote by $S^p$ the infinite-dimensional Schatten class of
exponent $p$. The definition and properties of the vector-valued
version of $S^p$ are similar to the finite-dimensional case, see
Chapter $1$ of \cite{P2}. We omit the proof of the following
result since the arguments to be used can be found in the proof of
Theorem \ref{BochnerLebesgue}.

\begin{The} \label{PisierSchatten} Let $1 \le p, q \le 2$ and
assume that $E$ is an operator space having Fourier type $p$ and
Fourier cotype $q'$ with respect to a compact group $G$. Then
\begin{itemize}
\item[$1.$ ] $S^r(E)$ has Fourier type $p$ with respect
to $G$ for all $p \le r \le p'$.
\item[$2.$ ] $S^s(E)$ has Fourier cotype $q'$ with
respect to $G$ for all $q \le s \le q'$.
\end{itemize}
Moreover, $\Co_p^1(S^r(E),G) = \Co_p^1(E,G)$ and
$\Co_{q'}^2(S^s(E),G) = \Co_{q'}^2(E,G)$.
\end{The}

\begin{Rem}
\emph{We already know that the Fourier type and cotype become
stronger conditions on the pair $(E,G)$ as $p$ and $p'$ approach
$2$. This gives rise to the notions of sharp Fourier type and
cotype exponents. The problem of finding the sharp exponents of a
given operator space is highly non-trivial even for the simplest
case of Lebesgue spaces or Schatten classes. Part of this problem
is solved in \cite{GMP}. Namely, if $1 \le p \le 2$ and $(\Omega,
\mathcal{M}, \nu)$ is not the union of finitely many $\nu$-atoms,
then we show that $L^p(\Omega)$ has sharp Fourier type $p$ with
respect to any compact semisimple Lie group. By duality we also
get that $L^{p'}(\Omega)$ has sharp Fourier cotype $p'$ for those
groups. By the nature of $\Omega$ and Proposition \ref{Subspace}
we have $$\mathcal{C}_q^1(L^p(\Omega),G) \ge \lim_{n \rightarrow
\infty} \,\ \mathcal{C}_q^1(l^p(n),G)$$ for $1 \le p < q \le 2$.
Moreover, Theorem \ref{Local} gives $\mathcal{C}_q^1(l^p(n),G) \le
n^{1/p - 1/q}$. The main result of \cite{GMP} asserts that there
exists a positive constant $\mathcal{K}(G,q)$, such that
$\mathcal{K}(G,q) \,\ n^{1/p - 1/q} \le \mathcal{C}_q^1(l^p(n),G)
\le n^{1/p - 1/q}$ for all $n \ge 1$ and any compact semisimple
$G$. The constant $\mathcal{K}(G,q)$ can be defined as
$$\mathcal{K}(G,q) = \inf_{n \ge 1} \sup \left\{
\frac{\|\widehat{f}\|_{\mathcal{L}^{q'}(\G)}}{\|f\|_{L^q(G)}}: \,\
f \ \ \mbox{central}, \ \ f \in L^q(G), \ \ \mbox{supp}(f) \subset
\mathcal{U}_n \right\}$$ where $\{\mathcal{U}_n: n \ge 1\}$
denotes a neighborhood basis at the identity of $G$. The
interesting point lies in the inequality $\mathcal{K}(G,q) > 0$
which constitutes a local variant of the Hausdorff-Young
inequality on $G$ with parameter $q$. The proof obtained for this
local inequality is based upon the semisimplicity of $G$ since it
uses the very well-developed theory of representations on such
kind of groups. The need to use these algebraic techniques forced
us to present the proof of this result in a separate work, see
\cite{GMP}.}
\end{Rem}

\bibliographystyle{amsplain}


\begin{flushleft}
\textsf{Jos\'{e} Garc\'{\i}a-Cuerva and Javier Parcet \\ Departamento
de Matem\'{a}ticas, C-XV, \\ Universidad Aut\'{o}noma de Madrid, \\
28049 Madrid,  Spain} \\ ${}$ \\ \texttt{jose.garcia-cuerva@uam.es} \\
\texttt{javier.parcet@uam.es}
\end{flushleft}

\end{document}